\documentclass[a4paper, 12pt]{article}
\usepackage[utf8]{inputenc}
\usepackage{a4}%kleinere Ränder
\usepackage[T1]{fontenc}%Silbentrennung
\usepackage[english]{babel}
\usepackage[hang,flushmargin]{footmisc} %No indentation at footnotes
\usepackage{amssymb,amsmath,amsthm}
\usepackage{enumerate}
\usepackage{graphicx}
\usepackage{epsfig}
\usepackage{comment}
\usepackage{xspace}
\usepackage{float}
\usepackage{multirow}
\usepackage{extarrows}
\usepackage{mathtools}
\usepackage{marvosym}
\usepackage{color}
\usepackage{nomencl}%Notationsverzeichnis
\usepackage{nicefrac}
\usepackage{cite}
\makenomenclature
\newcommand{\N}{\ensuremath{\mathbb{N}}\xspace}

\newcommand{\R}{\ensuremath{\mathbb{R}}\xspace}

\newcommand{\eps}{\epsilon}
\newcommand{\gr}[1]{(\ref{#1})}
\newcommand{\dx}{\,\mathrm{d}x}
\newcommand{\per}{\mathrm{Per}}
\newcommand{\leb}{\mathcal{L}}

\renewcommand{\epsilon}{\varepsilon}

\newcommand{\mres}{\mathbin{\vrule height 1.6ex depth 0pt width 0.13ex\vrule height 0.13ex depth 0pt width 1.3ex}} % 'restriction of a measure' sign
\renewcommand{\ell}{\mathfrak{l}}
\newcommand\blfootnote[1]{%
  \begingroup
  \renewcommand\thefootnote{}\footnote{#1}%
  \addtocounter{footnote}{-1}%
  \endgroup
}
\begin{document}

\numberwithin{equation}{section}
%Eigene Theoremumgebungen
\newtheoremstyle{break}{15pt}{15pt}{\itshape}{}{\bfseries}{}{\newline}{}
\theoremstyle{break}
\newtheorem*{Satz*}{Theorem}
\newtheorem*{Rem*}{Remark}
\newtheorem*{Lem*}{Lemma}
\newtheorem{Satz}{Theorem}[section]
\newtheorem{Rem}{Remark}[section]
\newtheorem{Exmp}{Example}[section]
\newtheorem{Lem}{Lemma}[section]
\newtheorem{Prop}{Proposition}[section]
\newtheorem{Cor}{Corollary}[section]
\theoremstyle{definition}
\newtheorem{Def}[Satz]{Definition}
\newtheorem*{Ass}{General Assumptions}
\parindent2ex

\begin{center}{\Large \bf  A reciprocity principle for constrained isoperimetric problems and existence of isoperimetric subregions in convex sets}
\end{center}

\begin{center}
Michael Bildhauer, Martin Fuchs, Jan M\"uller
\end{center}

\begin{abstract}
It is a well known fact that in $\R^n$ a subset of minimal perimeter $L$ among all sets of a given volume  is also a set of maximal volume among all sets of the same perimeter $L$. This is called the reciprocity principle for isoperimetric problems. The aim of this note is to prove this relation in the case where the class of admissible sets is restricted to the subsets of some subregion $G\subsetneq\R^n$. Furthermore, we give a characterization of those (unbounded) convex subsets of $\R^2$ in which the isoperimetric problem has a solution. The perimeter that we consider is the one relative to $\R^n$.
\end{abstract}


\begin{section}{Introduction}
In its classical form, the isoperimetric problem asks for the maximal area which can be enclosed by a curve of given length. In modern mathematical terms, the task is to determine a measurable subset of  $\R^n$ which has maximal Lebesgue measure $A$ among all sets of a given perimeter $L\in [0,\infty)$. Assuming the existence of a solution, Steiner in the first half of the 19$^{th}$ century showed by means of elementary geometric arguments that in $\R^2$  the only possible candidate for a solution is the circle of perimeter $L$ (which has already been suspected since antiquity). However, the existence part turns out to require a more subtle reasoning. In nowadays mathematics, it is usually treated in the framework of the theory of convex sets (see \cite{La}),  of  Caccioppoli sets (i.e.  sets of finite perimeter, see \cite{Fu} for a compact introduction to the topic) or integral rectifiable currents (see, e.g., \cite{Mo}). While the first method relies on an a priori establishment of the convexity of a solution, the second approach is based on a reformulation of the problem: it is easy to see by a scaling argument that a set in $\R^n$ which maximizes volume to a given perimeter is also a set with minimal  perimeter among all sets of the respective volume. Some authors refer to this connection between the two problems as ``reciprocity'' (see, e.g. \cite{HC}, \cite{Gian}). The situation is different if we restrict our class of admissible sets to those which lie inside a proper subset $G\subsetneq\R^n$. Then the classical arguments via convexity and scaling may fail depending on the geometry of $G$ and it is no longer clear that a volume maximizing subset with given perimeter occurs as a solution of the reciprocal ``minimal-perimeter''-problem. The second type of problem, i.e. finding sets of minimal perimeter of a given volume inside a subset of $\R^n$ is well addressed in the literature, ranging from existence results (which is clear if $G$ is bounded, see \cite{BB}, Theorem 1.2.2) up to  results  concerning the regularity of the boundary and convexity  (see, e.g., \cite{Giu2}, \cite{Ta}, \cite{Ta2}, \cite{Ta3}, \cite{GMT} and \cite{SZ}). In \cite{Be2}, Besicovitch investigated volume maximal subsets under the additional assumption of convexity, thereby avoiding the difficulties about the existence question. The authors, not being experts in this field, do not claim that all results concerning the reciprocity to be completely new, but as we did not find anything about it after an extensive literature research, we decided to give the proof here. 

 Actually, our considerations evolved from the following simple question (which is similar to the problem considered in \cite{Be1}): given the stripe $G=\R\times [0,1]$ in $\R^2$ and $L>\pi$, what is the shape of an area maximizing subset $E\subset G$ with perimeter $L$, or, if you prefer a more colloquial phrasing: what is the shape of the table with the largest surface area that fits inside a narrow room of rectangular layout under the condition that a given number of persons should be able to take a chair? (Obviously, in order to provide space for the chairs we then have to solve the problem in an inner parallel set of the room). It turns out that a solution exists in form of a rectangle with two semicircles of radius $1/2$ attached to two opposing sides (which might be the reader's intuitive guess). 
 
\section{Notation and statement of the results}
By a \textit{subregion of } $\R^n$, we mean a subset $G\subset\R^n$ ($n\geq 2$) such that 
\begin{align}\label{G}
G\;\text{ is open }\quad\text{ and }\quad  \leb^n(\partial G)=0 \tag{$\ast$}
\end{align}
($\leb^n$ denoting Lebesgue's measure), i.e. an open continuity set of Lebesgue's measure. For an arbitrary $\mathcal{L}^n$-measurable set $E$ and an open set $F\subset\R^n$ we define the perimeter of $E$ in $F$ by
\[
\per(E;F):=\sup\left\{\int_{E}\mathrm{div}\varphi \dx\,:\,\varphi\in C^1_0(F,\R^n),\,|\varphi|\leq 1\right\}.
\]
If $F=\R^n$, then we write  for short $\per(E)$ instead of $\per(E;\R^n)$.
\begin{Rem}
 Note that $\per(E)$ is allowed to take the value $+\infty$. If $\per(E)$ is finite, then it is well known (cf. \cite{AFP}) that the characteristic function $\chi_E$, viewed as an element of $L^1(\R^n)$, is a function of bounded variation, i.e. $\chi_E$ has a distributional derivative $D\chi_E$ in form of a Radon measure of finite total mass. If $\per(E;K)$ is finite for any open subset $K\Subset\R^n$, then we say that $E$ has locally finite perimeter, or that $E$ is a \textit{Caccioppoli set} (cf. \cite{Giu}). Note that if two sets $E$ and $F$ differ only by a Lebesgue null set, then their perimeters coincide.
\end{Rem}

We consider the following pair of problems:
\begin{align}
\tag{$P$}
\left\{\begin{aligned}
&\text{Given } L\in [0,\per(G)), \text{ find a subset }E\subset \overline{G}\\
&\text{ with }\per(E)=L\text{ and such that }\mathcal{L}^n(E)\text{ is maximal}
\end{aligned}\right\}
\end{align}
and the corresponding ``reciprocal'' problem
\begin{align}
\tag{$P_*$}
\left\{\begin{aligned}
&\text{Given } A\in [0,\mathcal{L}^n(G)), \text{ find a subset }E\subset \overline{G}\\
&\text{with }\mathcal{L}^n(E)=A\text{ and such that }\per(E)\text{ is minimal}.
\end{aligned}\right\}
\end{align}
Then we have the following result:
\begin{Satz}[Reciprocity]\label{reci}
Let $G\subset\R^n$ be a subregion in the sense of \gr{G}, $L\in [0,\per(G))$ and $A\in [0,\mathcal{L}^n(G))$.
\begin{enumerate}[i)]
\item If $E\subset \overline{G}$ is a solution of problem $(P)$, then $E$ has minimal perimeter among all subsets of $G$ which have the same volume $\leb^n(E)$.
\item If $F\subset \overline{G}$ is a solution of problem $(P_*)$, then $F$ has maximal volume among all subsets of $G$ which have the same perimeter $\per(F)$.
\end{enumerate}
\end{Satz}
Note that Theorem \ref{reci} does not say anything about the existence of a solution of problem $(P)$. What it does say, however, is that any solution of $(P)$ occurs as a solution of $(P_*)$ and vice versa. If we additionally assume $G$ to be bounded, the existence of a solution of problem $(P_*)$ follows easily: let $(F_k)_{k\in\N}$ be a perimeter minimizing sequence of subsets in $G$ with $\leb^n(F_k)=A$. Then the corresponding sequence of characteristic functions $(\chi_{F_k})$ is bounded in some space $BV(B_R(0))$, where $R>0$ is large enough s.t. $G\Subset B_R(0)$. By the BV-compactness Theorem (cf Theorem 3.23 in \cite{AFP}), there is a function $f\in BV(B_R(0))$ such that (at least for a subsequence)
\[
 \chi_{F_k}\rightarrow f \text{ in } L^1(\R^n)\quad\text{and}\quad |Df|(B_R(0))\leq \liminf\limits_{k\rightarrow\infty}\per(F_k).
\]
Since (after possibly passing to another subsequence) $\chi_{F_k}\to f$ a.e., $f$ is (up to a set of measure zero) the characteristic function of the set $E=\{x\in \overline{G}\,:\,f(x)=1\}$, which therefore is a solution of problem $(P_*)$. Of course this ``direct method'' fails if we apply it to a volume maximizing sequence $(E_k)_{k\in N}$ of subsets with fixed perimeter $\per(E_k)=L$, since the limit set $E\subset\overline{G}$ will in general  not satisfy $\per(E)=L$. Instead, we can use the reciprocity principle to show the existence of a solution of problem $(P)$:
\begin{Satz}[Existence in bounded subregions]\label{exi}
Let $G\subset\R^n$ be a subregion in the sense of \gr{G} which is additionally bounded. Then it holds:
\begin{enumerate}[i)]
\item For all $L\in [0,\per(G))$ problem $(P)$ admits a solution.
\item If $L>\per(G)$, then problem $(P)$ does not have a solution.
\end{enumerate}
\end{Satz}
\begin{Rem}
Note that in the case ``$L=\per(G)<\infty$'', the set $G$ itself is a trivial solution of $(P)$. 
\end{Rem}
Our last result concerns the existence of solutions of problem $(P)$ in unbounded subregions. In general, isoperimetric sets do not necessarily exist for arbitrary choices of $L\in [0,\per(G))$ as the following example shows:
\begin{figure}[H]
\begin{center}
\includegraphics[scale=1]{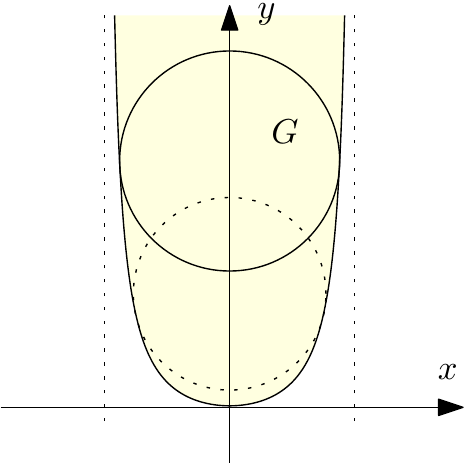}
\end{center}
\end{figure}
\begin{Exmp}\label{exmp} Consider $G:=\bigg\{(x,y)\in\R^2\,:\,-1<x<1,\;y\geq \frac{x^2}{1-x^2}\bigg\}$ (see the picture above) and choose $L=2\pi$. Then $(P)$ does not have a solution in $G$: by comparison with the sequence of maximal disks  with center $(0,n)$ inside $G$ ($n\in\N$), we see that the area of a solution of $(P)$ would be $\pi$. 
But by the isoperimetric inequality in $\R^2$, this can only be attained by a circle of radius $1$, which does not exist in $G$.
\end{Exmp}

However, in two dimensions we can give a complete characterization of the \textit{convex} subregions in which $(P)$ can be solved for any choice of $L\in[0,\per(G))$:

\begin{Satz}[Existence in convex subregions]\label{exconv}
Let $G\subset\R^2$ be open and convex and define
\[
\mathfrak{r}(G):=\sup\big\{r\in [0,\infty)\,:\,\exists x\in G\text{ such that }B_r(x)\subset\overline{G}\big\}.
\]
Then, if either
\begin{enumerate}[\hspace{1cm}(i)]
\item $\mathfrak{r}(G)=\infty\quad\text{ or }$
\item $\mathfrak{r}(G)=\max \big\{r\in [0,\infty)\,:\,\exists x\in G\text{ such that }B_r(x)\subset\overline{G}\big\}$,
\end{enumerate}
problem $(P)$ admits a solution for any $L\in [0,\per(G))$, which in addition is convex.
\end{Satz}

\begin{Rem}
If $G$ is a bounded convex subregion, then the existence of a maximal disk in $\overline{G}$, i.e. a disk of radius $\mathfrak{r}(G)$ (also called \textit{incircle}, or \textit{inball} in dimensions $n\geq 3$) follows from general principles, cf.  Corollary 16.2 in \cite{La}.
\end{Rem}

\begin{Rem}
If we additionally demand $G$ to satisfy the so called ``great circle condition'' from \cite{SZ} (i.e. there exists an inball $B\subset G$ and a hyperplane $H\subset\R^n$, containing the center of $B$, such that $\partial B\cap H\subset \partial G$), then Theorem 3.31 in \cite{SZ} together with our Lemma \ref{Lemtr} suffices to establish the assertion of Theorem \ref{exconv} in arbitrary dimensions $n\geq 3$ (cf. also Remark \ref{Rem5.1} at the end of the proof of Theorem \ref{exconv}). It is an open question if one can drop the great circle condition in this context.
\end{Rem}

\begin{Rem}
We would like to emphasize at this point that in contrast to e.g. \cite{LRV} we consider the perimeter of sets relative to $\R^n$.
\end{Rem}

\begin{Rem}\label{Rem2.5}
\begin{enumerate}[i)]
\item The existence of a solution in the case $\mathfrak{r}(G)=\infty$ in form of a disk is clear.
\item As Example \ref{exmp} shows, the theorem above gives  a sharp characterization of the convex sets in $\R^2$ in which $(P)$ (and thus, by Theorem \ref{reci} $(P_\ast)$) can be solved.
\item Theorem \ref{exconv} particularly applies to our example of a ``maximal table in a narrow room'' from the introduction. By Theorem 3.31 in \cite{SZ}, we further obtain that the solution of $(P_\ast)$ in the stripe $S=\R\times [0,1]$ is the convex hull of two incircles, i.e. has the shape of a ``stadium'' (compare also *1.4.3 on p. 5 in \cite{BSZ}):  
\begin{figure}[H]
\begin{center}
\includegraphics{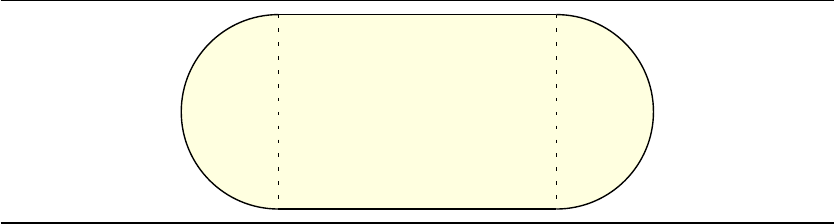}
\end{center}
\end{figure}
\end{enumerate}
\end{Rem} 
\end{section}

\section{Proof of Theorem \ref{reci}}
We start with part i). Let $L\in [0,\per(G))$ be given and let $E\subset\overline{G}$ have maximal volume $A=\leb^n(E)$ among all subsets of $\overline{G}$ with perimeter $L$. Assume that there is another subset $E'\subset\overline{G}$ with $\leb^n(E')=A$ and $\per(E')=L'<L=\per(E)$. Due to $\per(E')<\per(G)$, it must hold
\[
 \leb^n(G-E')>0.
\]
Thus, by Lebesgue's density Theorem (see \cite{EG2}, Corollary 3 on p. 45) there exists a point $x_0\in G-E'$ for which
\[
 \lim_{\rho\downarrow 0}\frac{\leb^n\big(B_\rho(x_0)\cap (G-E')\big)}{\leb^n(B_\rho(x_0))}=1
\]
and therefore
\[
 \lim_{\rho\downarrow 0}\frac{\leb^n\big(B_\rho(x_0)\cap E'\big)}{\leb^n(B_\rho(x_0))}=0.
\]
Now choose $\rho_0>0$ so small that
\[
 \per(B_{\rho_0}(x_0))<\frac{L-L'}{2}\quad\text{ and }\quad \leb^n(B_{\rho_0}(x_0)\cap E')<\frac{1}{2}\leb^n(B_\rho(x_0))
\]
and consider the set
\[
 \widetilde{E}:=E'-B_{\rho_0}(x_0).
\]
Then it holds
\[ 
\per(\widetilde{E})\leq \per(E')+\per(B_{\rho_0}(x_0))<L'+\frac{L-L'}{2}<L.
\]
Set $\widetilde{L}:=L-\per(\widetilde{E})$ and note that
\[
 \widetilde{L}>\frac{L-L'}{2}>\per(B_{\rho_0}(x_0)).
\]
Therefore, we can choose a compact subset $C\Subset B_{\rho_0}(x_0)$ such that
\[
 \leb^n(C)>\frac{1}{2}\leb^n(B_{\rho_0}(x_0))\quad\text{ and }\quad \per(C)=\widetilde{L}.
\]
But then the set $\hat{E}:=\widetilde{E}\cup C$ satisfies
\[
 \per(\hat{E})=\per(\widetilde{E})+\per(C)=L=\per(E)
\]
as well as
\[
 \leb^n(\hat{E})=\leb^n(\widetilde{E})+\leb^n(C)>\leb^n(E),
\]
which contradicts the volume maximality of the subset $E$.

We proceed with the proof of part ii). Let $A\in [0,\leb^n(G))$ be given and let $F\subset\overline{G}$ have minimal perimeter $L=\per(F)$ among all subsets of $\overline{G}$ with volume $A=\leb^n(F)$. Assume that there is another subset $F'$ with $\per(F')=L$ and $\leb^n(F')=A'>A$. For $\alpha\in \R$ consider the half space
\[
 H_\alpha:=\big\{x\in\R^n\,:\, x\cdot e_n\geq \alpha\big\},
\]
where $e_n:=(0,...,0,1)^T$ and define the function
\[
 h:\R\rightarrow\R,\,\alpha\mapsto \leb^n(F'\cap H_\alpha).
\]
Since by Lebesgue's theorem on dominated convergence
\[
h(\alpha)=\int_{F'}\chi_{H_\alpha}\dx=\lim_{n\to\infty}\int_{F'}\chi_{H_{\alpha_n}}\dx\quad\text{for any sequence }\alpha_n\to\alpha,
\]
the function $h$ is continuous and 
\[
\lim\limits_{\alpha\to-\infty}h(\alpha)=A'\quad\text{ and }\quad \lim\limits_{\alpha\to\infty}h(\alpha)=0.
\]
Thus, by the intermediate value Theorem, there exists $\alpha_0\in \R$ with $h(\alpha_0)=A$. 
\begin{Lem}\label{Lemh}
 With $\alpha_0$ and $H_\alpha$ as above, it holds
\[
 \per(F'\cap H_{\alpha_0})<\per(F')=L.
\]
\end{Lem}

Thus, we see that the set $F'\cap H_{\alpha_0}$ contradicts the perimeter minimality of the set $F$. It remains to give a proof of the above assertion.

\noindent\textit{Proof of the lemma}. Without loss of generality we may assume $\alpha_0=0$. Set
\begin{gather*}
 H_+:=\big\{x\in\R^n\,:\,x\cdot e_n>0\big\},\\
 H_0:=\big\{x\in\R^n\,:\,x\cdot e_n=0\big\},\\
 H_-:=\big\{x\in\R^n\,:\,x\cdot e_n<0\big\}.
\end{gather*}
Due to $\leb^n(F'\cap H_0)=0$, we have $\per(F')=\per\big((F'\cap H_+)\cup (F'\cap H_-)\big)$. Let further $\chi_+,\chi_-\in L^1(H_0)$  be the trace of $\chi_{F'\cap H_+}$ and $\chi_{F'\cap H_-}$ in $H_0$, respectively, in the sense of $BV$-functions. (Note that these exist by Theorem 3.77 in \cite{AFP} and are summable w.r.t. $\mathcal{H}^{n-1}\mres H_0$ since $F'\cap H_+$ and $F'\cap H_-$ have finite perimeter in $\R^n$). Then it holds
\begin{align*}
 \per(F')=\per(F';H_+)+\per(F';H_-)+\int_{H_0}|\chi_+-\chi_-|\,\mathrm{d}\mathcal{H}^{n-1}
\end{align*}
as well as
\begin{align*}
 \per(F'\cap H_+)=\per(F';H_+)+\int_{H_0}|\chi_+|\,\mathrm{d}\mathcal{H}^{n-1}.
\end{align*}
It therefore follows
\begin{align*}
 \per(F'\cap H_+)\leq \per(F';H_+)+\per(F';H_-)+\int_{H_0}|\chi_+-\chi_-|\,\mathrm{d}\mathcal{H}^{n-1}\\
+\int_{H_0}|\chi_-|\,\mathrm{d}\mathcal{H}^{n-1}-\per(F';H_-)\\
=\per(F')-\left[\per(F';H_-)-\int_{H_0}|\chi_-|\,\mathrm{d}\mathcal{H}^{n-1}\right].
\end{align*}
It thus remains to show that $\per(F';H_-)-\int_{H_0}|\chi_-|\,\mathrm{d}\mathcal{H}^{n-1}>0$. Note that due to our choice of $\alpha_0$ in Lemma \ref{Lemh} it holds
\[
 \leb^n(F'\cap H_-)=A'-A>0.
\]
In \cite{Mag}, Theorem 19.15 it is shown that the  following constrained isoperimetric problem in the half space,
\begin{align*}
 \text{``minimize }\per(S;H_-)\text{ among all sets }S\subset H_-\text{ with }\leb^n(S)=v\\
\text{ and }\per(S;H_0)=\sigma"
\end{align*}
for some given $v>0$, $\sigma\geq 0$ is solved by the segment of a ball $B_R\big((0,...,0,x_n)\big)$, where $R,x_n\geq 0$ are such that $\leb^n(B_R\cap H_-)=v$ and $\mathcal{H}^{n-1}(B_R\cap H_0)=\sigma$. But for a ball which intersects $H_-$ in a set of positive measure it is surely true that $\per(B_R;H_-)-\per(B_R;H_0)>0$. Hence, choosing 
\[
v=A'-A,\;\sigma=\int_{H_0}|\chi_-|\,\mathrm{d}\mathcal{H}^{n-1}
\]
 we infer that $\per(F';H_-)-\int_{H_0}|\chi_-|\,\mathrm{d}\mathcal{H}^{n-1}>0$. \qed

\section{Proof of Theorem \ref{exi}}
i). We start with the observation that the so called \textit{isoperimetric profile} of $G$, i.e. the function
\[
 \mathfrak{l}:[0,\leb^n(G))\rightarrow \R,\, \mathfrak{l}(A):=\inf\big\{\per(E)\,:\,E\subset\overline{G}\text{ with }\leb^n(E)=A\big\}
\]
is well-defined for any open subset $G\subset\R^n$. Our further proof relies on the following properties of the function $ \mathfrak{l}$ on bounded subregions, which might be well known in the more general context of compact Riemannian manifolds (see, e.g., \cite{Ri}):
\begin{Lem}\label{Leml}
 Let $G\subset\R^n$ be a bounded subregion and let $ \mathfrak{l}:[0,\leb^n(G))\rightarrow \R$ be defined as above. Then it holds
\begin{enumerate}[a)]
 \item The function $ \mathfrak{l}$ is strictly increasing,
 \item $ \mathfrak{l}$ is continuous,
 \item $ \mathfrak{l}$ maps the interval $[0,\leb^n(G))$ bijectively to the interval $[0,\per(G))$.
\end{enumerate}
\end{Lem}

Before we prove the lemma, let us see how it applies to the proof of Theorem \ref{exi}. Assume that $L\in [0,\per(G))$ is given. Then, by part c) of the Lemma there is $A\in [0,\leb^n(G)$ such that $ \mathfrak{l}(A)=L$. Let $E\subset\overline{G}$ be a subset with $\leb^n(E)=A$ and $\per(E)= \mathfrak{l}(A)=L$ (which exists due to the boundedness of $G$). Then $E$ has maximal volume among all subsets with perimeter $L$, because the existence of any set with strictly larger measure and the same perimeter $L$ would be in contradiction to the monotonicity of $ \mathfrak{l}$ from part a). Hence $E$ is a solution of problem $(P)$.

\noindent\textit{Proof of the lemma}. a). Let $0\leq A<A'<\leb^n(G)$. Since $G$ is bounded, we find subsets $E$ and $E'$ of $\overline{G}$ of volume $A$ and $A'$, such that $\per(E)=\ell(A)$ and $\per(E')=\ell(A')$. Then, as in Lemma \ref{Lemh} from the previous section, we can choose a real number $\alpha_0$ such that the intersection of $E'$ with the half space $H_{\alpha_0}$ has volume $A$. But, as it was shown in Lemma \ref{Lemh}, it then follows
\[
 \ell(A)\leq \per(E'\cap H_\alpha)<\per(E')=\ell(A').
\]

\noindent b). Let $A\in [0,\leb^n(G))$ be arbitrary and let $(A_k)_{k\in\N}$ be a sequence in $[0,\leb^n(G))$ which converges to $A$ from below. From the monotonicity of $\ell$, we infer that the limit of $\ell(A_k)$ exists and satisfies $\lim\limits_{k\rightarrow\infty}\ell(A_k)\leq \ell(A)$. Let $E_k\subset\overline{G}$ be such that $\per(E_k)=\ell(A_k)$. By the BV-compactness property (Theorem 3.23 in \cite{AFP}) there is a subset $E\subset\overline{G}$ such that (at least for a subsequence) 
\[
\chi_{E_k}\to\chi_E\text{ in }L^1(\overline{G})\quad\text{ and }\quad\per(E)\leq \liminf_{k\rightarrow\infty}\per(E_k)=\lim_{k\rightarrow\infty}\per(E_k).
\]
Now if $\lim\limits_{k\rightarrow\infty}\ell(A_k)< \ell(A)$, then the set $E$ would satisfy $\leb^n(E)=A$ and $\per(E)<\ell(A)$, which contradicts the definition of $\ell(A)$.

Let now $(A_k)$ converge to $A$ from above and let $\eps>0$ be given. Let $E\subset\overline{G}$ be such that $\leb^n(E)=A$ and $\per(E)=\ell(A)$. Since $\leb^n(G-E)>0$, Lebesgue's density theorem implies the existence of a point $x_0\in G-E$ for which
\[
 \lim_{\rho\downarrow 0}\frac{\leb^n(B_\rho(x_0)\cap E)}{\leb^n(B_\rho(x_0))}=0.
\]
Choose $\rho_0>0$ so small that $\leb^n(B_\rho(x_0)\cap E)<\frac{1}{2}\leb^n(B_\rho(x_0))$ and such that $\per(B_{\rho_0}(x_0))<\eps$. Consider the set
\[
 E':=E\cup B_{\rho_0}(x_0).
\]
Then it holds
\[
 \leb^n(E'):=A'>A\quad\text{ and }\quad \per(E')\leq \per(E)+\per(B_{\rho_0}(x_0))\leq \ell(A)+\eps.
\]
Thus, if we choose $N\in \N$ large enough such that $A_k<A'$ for all $k>N$, it follows from the monotonicity of $\ell$ that
\[
 \ell(A)\leq \ell(A_k)\leq \ell(A')\leq \ell(A)+\eps,
\]
and b) of Lemma \ref{Leml} is proved.

\noindent c). Having established parts a) and b), it suffices to prove $\ell(A)\rightarrow \per(G)$ as $A\rightarrow \leb^n(G)$. But this follows easily from the lower semicontinuity of the perimeter with respect to $L^1$-convergence: let $A_k\to \leb^n(G)$ in $[0,\leb^n(G))$ and choose $E_k$ with $\per(E_k)=\ell(A_k)$. Then $\chi_{E_k}\to \chi_G$ in $L^1(\R^n)$ and therefore 
\[
\per(G)\leq \liminf\limits_{k\rightarrow\infty}\per(E_k)=\liminf\limits_{k\rightarrow\infty}\ell(A_k)\leq \per(G). 
\]
\qed
\newpage

We now come to the proof of ii). Let $L>\per(G)$ be given. Assume that $E\subset \overline{G}$ has maximal volume among all subsets of perimeter $L$. Let $x_0\in G$ be some point and choose $r>0$ such that $B_r(x_0)\subset G$. For any $k\in\N$, let $C_k\Subset B_{r/k}(x_0)$ denote a compact subset with perimeter $\per(C_k)=L-\per(G)$ and consider the sequence
\[
E_k:=G-C_k
\] 
of subsets of $G$. Then $\per(E_k)=\per(G)+\per(C_k)=L$ and $\leb^n(G-E_k)\rightarrow 0$ for $k\to \infty$. Therefore it must hold $\leb^n(E)=\leb^n(G)$ and thus $\leb^n(\overline{G}-E)=0$. But then $\per(E)=\per(\overline{G})=\per(G)$ in contradiction to our assumption that $E$ solves $(P)$ for $L>\per(G)$. \qed

\section{Proof of Theorem \ref{exconv}}

Let $G\subset\R^2$ be an open convex set and let w.l.o.g. $B_{1}(0)\subset \overline{G}$ be a maximal disk in $\overline{G}$. We may assume $L>2\pi$ (since otherwise the disk $B_{L/(2\pi)}(x)\subset B_1(0)$ would be a trivial solution of problem $(P)$) as well as $\per(G)=\infty$ (since otherwise $\mathrm{diam}(G)\leq \frac{1}{2}\per(G)<\infty$ and existence in bounded subregions is clear by Theorem \ref{exi}).  For $k\in\N$ we define 
\[
 G_k:=\overline{G}\cap [-k,k]\times [-k,k],
\]
which is a convex and bounded subset of $\R^2$. Let
\[
 A:=\sup\big\{\leb^2(E)\,:\, E\subset \overline{G},\;\per(E)=L\big\}\leq \frac{L^2}{4\pi}.
\]
We claim the following:
\begin{Lem}\label{Lemtr}
 Choose $k\in \N$ large enough such that $\per(G_k)>L$ and let $E_k\subset G_k$ be sets of maximal area among all subsets of $G_k$ with perimeter $L$ (note that such a set exists by Theorem \ref{exi} i)). Then 
\[
 \lim_{k\rightarrow\infty}\leb^2(E_k)=A,
\]
i.e. $(E_k)$ is an area-maximizing sequence in $G$.
\end{Lem}

\noindent\textit{Proof of the lemma}. Let $\eps>0$ be given. Choose a subset $F\subset \overline{G}$ of perimeter $L$, for which $\leb^2(F)>A-\frac{\eps}{2}$ and $k_0\in \N$ large enough, such that 
\[
\leb^2(F)-\leb^2\big(F\cap [-k_0,k_0]\times[-k_0,k_0]\big)<\frac{\eps}{2}.
\]
 We set $F_{k_0}:=F\cap [-k_0,k_0]\times[-k_0,k_0]$. Then, since intersecting $F$ with a square is the same as intersecting $F$ gradually with four half spaces, Lemma \ref{Lemh} yields
\[
 \per(F_{k_0}):=L'\leq \per(F)=L.
\]
From the proof of Reciprocity Theorem \ref{reci}, we see that on the bounded subregion $G_k$ the function 
\[
 \mathfrak{a}_k:[0,\per(G_k))\to \R,\; L\mapsto \sup\big\{\leb^2(E)\,:\,E\subset G_k,\; \per(E)=L\big\}
\]
is the inverse of the corresponding function
\[
 \ell_k:[0,\leb^2(G_k))\to \R,\; A\mapsto \inf\big\{\per(E)\,:\,E\subset G_k,\; \leb^2(E)=A\big\}
\]
and thus, by Lemma \ref{Leml}, it is strictly increasing. Therefore we have
\[
 \leb^2(E_{k_0})=\mathfrak{a}_{k_0}(L)\geq \mathfrak{a}_{k_0}(L')\geq \leb^2(F_{k_0})\geq A-\eps.
\]
The result now follows since the sequence $\leb^2(E_k)$ is increasing. \qed

\begin{Lem}\label{Lem5.2}
The sets $E_k$ from Lemma \ref{Lemtr} are convex. In particular, $\mathrm{diam}(E_k)\leq \frac{L}{2}$ for all $k\in\N$. 
\end{Lem}

\noindent\textit{Proof of the lemma}. By Theorem 1.1, each set $E_k$ is also a solution of the corresponding reciprocal problem in the convex and bounded set $G_k$. Convexity of the $E_k$ thus follows from quoting Theorem 3.24 in \cite{SZ} (see also Remark 3.25). For a plane convex set $E$, the inequality $\mathrm{diam}(E)\leq \frac{\per(E)}{2}$ is a triviality. (Note that this is wrong in $\R^n$ for $n\geq 3$). \qed

The further idea of the proof  is the following: in a convex set $G$, which contains a maximal disk, i.e. a disk of radius $\mathfrak{r}(G)$, any bounded convex subset $E\subset G$ has a translate $E'$ (i.e. $E'=E+p$ for some $p\in\R^2$) such that $E'$ intersects the maximal disk (see Lemma \ref{conv} below). Hence we may assume that the sets $E_k$ from Lemma \ref{Lemtr} all lie within some set $G_{K}$ for a fixed integer $K>L+1$. Let now $E\subset G_K$ be a set of maximal area among all subsets of $G_K$ which have perimeter $L$. Then $\leb^2(E)\geq \leb^2(E_k)$ for all $k\in\N$ and by Lemma \ref{Lemtr}, the constant sequence $(E)$ is seen to be an area-maximizing sequence of subsets in $G$. Hence $E$ is a solution of the problem $(P)$ in $G$. It remains to prove the following result on convex sets:

\begin{Lem}\label{conv}
Let $K\subset\R^2$ be a closed convex set and $B\subset K$ any maximal disk in $K$. Then every convex subset $C\subset K$ has a translate in $K$ which intersects $B$, i.e. there exists $p\in \R^2$ such that
\[
p+C\subset K\quad\text{ and }\quad (p+C)\cap B\neq \emptyset.
\]
\end{Lem}

\noindent\textit{Proof of the lemma}. We may assume that $B=\overline{B_1(0)}$ and that $C\cap B=\emptyset$. Let $H\subset\R^2$ be a common line of support of both $B$ and $C$, i.e. $H$ is a line which intersects $B$ as well as $\overline{C}$ and such that $B,C$ lie in the same of the two halfspaces which are separated by $H$ (see the picture below).

\begin{figure}[H]
\begin{center}
\includegraphics[scale=0.6]{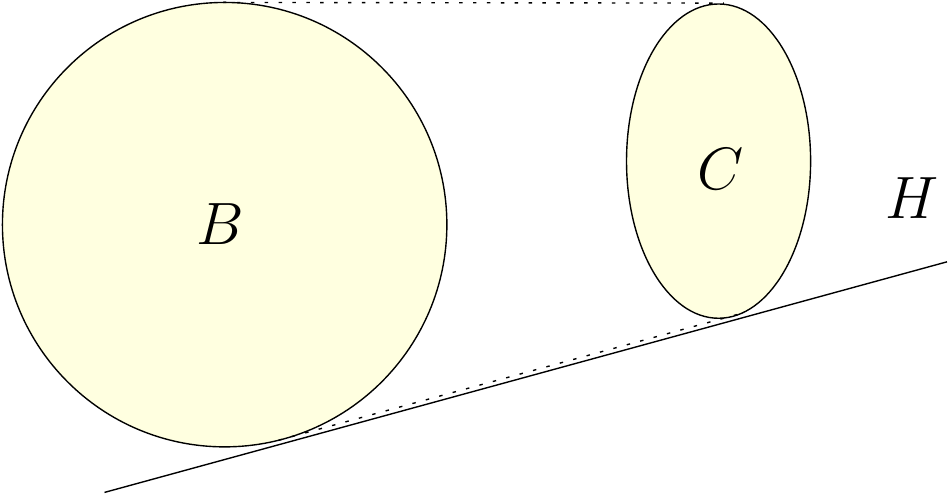}
\end{center}
\end{figure}

Such a line $H$ can always be found by taking a line of support at any boundary point of the convex hull $\mathrm{conv}(B\cup C)$, which neither lies on $\partial B$ nor on $\partial C$. W.l.o.g. we may assume that $H$ is the line $y=-1.$
Then we claim that $C$ lies completely within the stripe
\[
S:=\big\{(x,y)\in \R^2\,:\,|y|\leq 1\big\}.
\]
Indeed, let $a\in H\cap\overline{C}$ and assume that there exists a point $b\in C-S$. Note that due to the convexity of $C$ and our assumption $B\cap C=\emptyset$, the line $\overline{ab}$ does not intersect $B$. Let $K'$ be the convex hull of the disk $B$ and the points $a,b$. Due to convexity, we have $K'\subset$ $K$. Furthermore, we see that $B$ cannot be a maximal disk in $K'$: since $\eps:=\mathrm{dist}(B,\overline{ab})>0$, we can translate $B$ by $\frac{\eps}{2}$ in $x$ direction, so that the shifted disk $B'$ intersects the boundary of $K'$ only in the point $(\eps/2,-1)$.  Due to $b\notin S$, $B'$ has positive distance to $\partial K'-H$ and can therefore not be maximal in $K'$ (see the picture below), which is in contradiction to our assumption.
\begin{figure}[H]
\begin{center}
\includegraphics[scale=0.6]{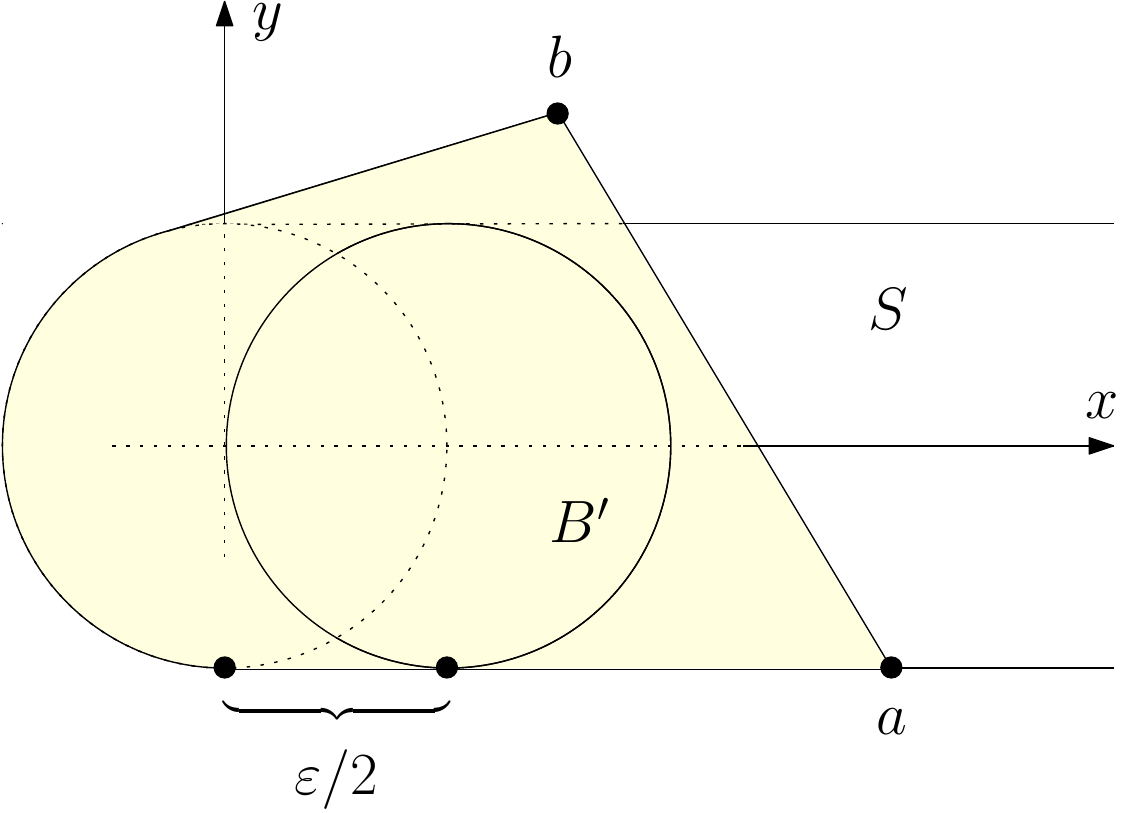}
\end{center}
\end{figure}
It follows that $C\subset S$. Let $d:=\mathrm{dist}(B,C)$. Then, by the convexity of $K$ the translate
\[
C-(d+1)\begin{pmatrix}1\\0\end{pmatrix}
\]
lies inside $K$ and intersects $B$, as claimed. \qed

\begin{Rem}\label{Rem5.1}
We would like to thank the referee for making us aware of the following fact: if the convex set $C$ in the proof of Lemma \ref{conv} is unbounded, then, by connecting points in $\{0\}\times [-1,1]$ with points $(x,y)$ for increasing values of $|x|$, one sees that either the half stripe $[0,\infty)\times [0,1]$ or $(-\infty,0]\times [0,1]$ must be contained in $C$. In particular, an isoperimetric subregion can be moved to lie inside this half stripe, which implies that it has the shape depicted in Remark \ref{Rem2.5} iii). Similar arguments show that an unbounded convex set in $\R^n$ for $n\geq 3$, which satisfies the ``great circle condition'' (cp. Remark \ref{Rem2.5}), coincides (up to a bounded part) with a spherical (half) cylinder, so that any isoperimetric subregion must be the convex hull of two inballs in this cylinder by Theorem 3.31 in \cite{SZ}. Note that this is false in general if the great circle condition is dropped (take, as an example, the set $[0,1]\times [0,1]\times \R$ in $\R^3$).
\end{Rem}

\begin{Rem}
 We would like to note that the method described above probably applies to a larger class of two-dimensional subregions $G$. We say that $G\subset\R^2$ has the ``\textit{translation property}'', if $G$ contains a maximal disk $B$ such that any connected subset $E\subset G$ has a translate $E'=E+p$ in $G$ which intersects $B$. An example of a non-convex set which has the translation property is shown in the following picture:
\begin{figure}[H]
\begin{center}
\includegraphics[scale=1]{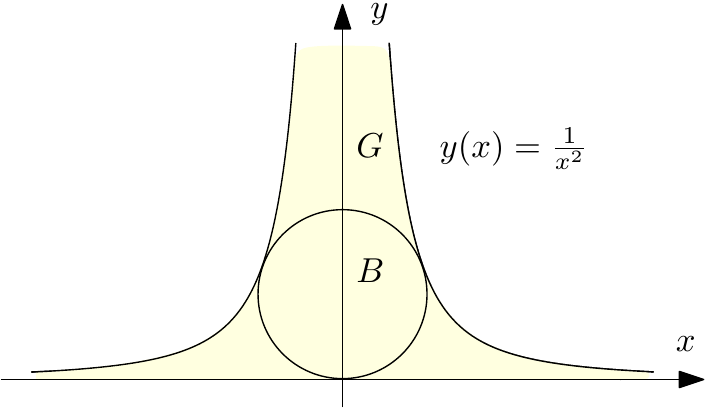}
\end{center}
\end{figure}
Then, if $G$ is a subregion which has the translation property and, in addition, $C^1$-smooth boundary, we can repeat the argument from the proof of Theorem \ref{exconv} for each connected component of an area maximizing sequence in $G$, which by \cite{Ta} is $C^1$-smooth as well. However, it seems to be difficult to give a sharp characterization of sets having the translation property in geometrical terms.
\end{Rem}

\blfootnote{
\begin{large}\Letter\end{large} Michael Bildhauer (bibi@math.uni-sb.de), Martin Fuchs (fuchs@math.uni-sb.de), Jan Müller (jmueller@math.uni-sb.de)\\
Saarland University (Department of Mathematics), P.O. Box 15 11 50, 66041 Saarbrücken, Germany }

\end{document}